\documentclass{amsart}
\usepackage{graphicx}

\newtheorem{thm}{Theorem}
\newtheorem{lem}{Lemma}

\begin{document}

\title{An asymptotic Robin inequality}
\author{Patrick Sol\'e}
\address{CNRS/LTCI, Universit\'e Paris-Saclay, 75 013 Paris, France.}
%\address{$^*$, King Abdulaziz University, Jeddah, Saudi Arabia.}
\email{sole@telecom-paristech.fr}
\author{Yuyang Zhu}
\address{Department of Math and Physics, Hefei University, Hefei, China.}
%\address{$^*$, King Abdulaziz University, Jeddah, Saudi Arabia.}
\email{zhuyy@hfuu.edu.cn}
\begin{abstract} The conjectured Robin inequality for an integer $n>7!$ is $\sigma(n)<e^\gamma n \log \log n,$ where $\gamma$ denotes Euler constant, and
$\sigma(n)=\sum_{d | n} d $. Robin proved that this conjecture is equivalent to Riemann hypothesis (RH).
Writing $D(n)=e^\gamma n \log \log n-\sigma(n),$ and $d(n)=\frac{D(n)}{n},$ we prove unconditionally that 
$\liminf_{n \rightarrow \infty} d(n)=0.$ 
The main ingredients of the proof are an estimate for Chebyshev summatory function,
and an effective version of Mertens third theorem due to Rosser and Schoenfeld. 
A new criterion for RH depending solely on $\liminf_{n \rightarrow \infty}D(n)$ is derived.
%We study the symmetry of cyclic combinatorial structures like cyclic codes, quasi-cyclic codes, Cayley graphs of cyclic groups. We explore the problem of equivalence testing, enumeration of non equivalent objects with the same parameters,  and of the existence of non trivial extra automorphisms. 
\end{abstract}

\maketitle
\let\thefootnote\relax
%%%%%%%%%%%%%%%%%%%%%%%%%%%%%%%%%%%%%%
\section{Introduction}
\subsection{History}
The conjectured Robin inequality for an integer $n>7!=5040$ is $\sigma(n)<e^\gamma n \log \log n,$ where $\gamma\approx 0.577\cdots$ denotes Euler constant, and $\sigma$ is the sum-of-divisors functions
$\sigma(n)=\sum_{d | n} d $. This inequality has been shown to hold unconditionally for families of integers that are
\begin{itemize}
 \item odd $>9$ \cite{C}
 \item square-free  $>30$ \cite{C}
 \item a sum of two squares and $>720$ \cite{BM}
 \item not divisible by the fifth power of a prime \cite{C}
 \item not divisible by the seventh power of a prime \cite{PS}
 \item not divisible by the eleventh power of a prime \cite{BT}
 
\end{itemize}
Ramanujan showed that Riemann Hypothesis implied that conjecture \cite{Ra}. Robin proved the converse statement \cite{R}, thus making
that conjecture a criterion for RH. This criterion was made popular by \cite{L} which derives an alternate criterion involving
Harmonic numbers.

\subsection{Contribution}
Denote the difference between the right hand side and the left hand side of Robin inequality by $D(n)=e^\gamma n \log \log n-\sigma(n).$ Let 
$d(n)=\frac{D(n)}{n}.$ The main result of this note is

{\thm For large $n$ we have $\liminf_{n \rightarrow \infty} d(n)=0.$}

Its proof will depend on the following intermediate result.

{\thm For large $n$ the quantity $\liminf_{n \rightarrow \infty} d(n)$ is finite and $ \ge 0.$}

The main ingredients of the proof of the latter are a combinatorial inequality between arithmetic functions (Lemma 1), 
an effective version of Mertens third theorem due to Rosser and Schoenfeld (Lemma 2), and an asymptotic estimate of
Chebyshev first summatory function (Lemma 4). Also needed is a result of Ramanujan of 1915 first published in 1997 \cite{Ra}.

We also study the asymptotic behavior of $D(n).$ Recall that a number is Colossally Abundant (CA) if it is a left-to-right maxima
for the function with domain the integers $x \mapsto \frac{\sigma(x)}{x^{1+\epsilon}},$ where $\epsilon$ is a real parameter.
Thus $n$ is CA iff $m<n \Rightarrow \frac{\sigma(m)}{m^{1+\epsilon}}<\frac{\sigma(n)}{n^{1+\epsilon}}.$
%Surprizingly the same statement as Theorem 2 for $D(n)$ is equivalent to RH.
{\thm We have the following limits when $n$ ranges over Colossally Abundant numbers.

\begin{itemize}
        \item If RH is false then $\liminf_{n \rightarrow \infty} D(n)=-\infty$ 
        \item If RH is true then $\lim_{n \rightarrow \infty} D(n)=\infty.$
       \end{itemize}

}

This result constitutes a new criterion for RH.
Its proof will depends, for the RH false part, on an oscillation theorem of Robin \cite{R}, modelled after and depending upon an oscillation theorem of Nicolas \cite{N}
for the Euler totient function. For the RH true case, we use a result of Ramanujan from 1915, first published in 1997 \cite{Ra}.
\subsection{Organization} The material is arranged as follows. The next section contains the proof of Theorem 1, Section 3 that of Theorem 2, and Section 4 that of Theorem 3. Section
5 concludes and gives some open problems.
\section{Proof of Theorem 1}
The result will follow from Theorem 2 if we exhibit a sequence of integers $n_m$ with $\lim_{m \rightarrow \infty} D(n_m)=0.$ We follow the approach
of \cite[\S 4, proof of Lemma 4.1, 1), p. 366]{C}. Consider $n$ of the shape $n=\prod_{p\le x} p^{t-1},$ with $t>1$ integer and $x$ real, both going infinity, and to be specified later. 
By this reference, we have
$$d(n)=e^\gamma  \log \log n (1-\frac{1}{\zeta(t)}+o_t(1)),$$
with $\zeta$ Riemann zeta function. The error term can be made effective as follows. By \cite[(3.28),(3.30)]{RS} we have
$$ e^\gamma \log x (1-\frac{1}{2\log^2 x})\le\prod_{p\le x} (1-\frac{1}{p})^{-1}\le e^\gamma \log x (1+\frac{1}{\log^2 x}). $$

From the Euler product of $\zeta$ and \cite[Lemma 6.4]{C} we derive

$$ \frac{1}{\zeta(t)}\le\prod_{p\le x} (1-\frac{1}{p^t})\le\frac{\exp(\frac{t x^{1-t}}{t-1})}{\zeta(t)}.$$

Combining these four bounds together we can take $o_t(1)=O\big(\exp(\frac{t x^{1-t}}{t-1})-1\big)=O(x^{1-t}) .$
Now it is elementary to show for $t$ integer $> 1$ that $\zeta(t)=1+h(t),$ with $h(t)=O(1/2^t).$ Indeed

$$\frac{1}{2^t} \le \zeta(t)-1\le \sum_{m=1}^\infty\frac{1}{2^{mt}}=\frac{2^{-t}}{1-2^{-t}}. $$
Thus, reporting, we get
$$d(n)=e^\gamma  \log \log n \big(O(1/2^t)+O(x^{1-t})\big).$$
To achieve $d(n) \rightarrow 0,$ we need both $\log \log n <<2^t,$ and $\log \log n <<x^{t-1}.$
This is ensured if we take $x=p_m,$ and $t=m+1.$ 
In that case we have $\log \log n =\log m+\log \theta(p_m).$ By Lemma 4 below, $\log \theta(p_m)\sim \log p_m.$
On the other hand, $p_m \sim m \log m$ as is well-known (see e.g. \cite{EM}).
Combining the last two estimates we see that 
$\log \log n\sim 2\log m << 2^m.$ Similarly,
$\log \log n<<p_m ^m.$
\section{Proof of Theorem 2}
If $\liminf_{n \rightarrow \infty}d(n)=\infty$ then $\lim_{n \rightarrow \infty}D(n)=\infty,$ and by Robin criterion RH holds.
We know then by \cite[p.25]{Ra} that the sequence $d(n)\sqrt{\log n}$ admits finite upper and lower limits when $n$ ranges over CA numbers (see \S 4), which is a contradiction.

Assume therefore that $\liminf_{n \rightarrow \infty}d(n)$ is finite and let us show that it is $\ge 0.$
For any integer $n$ write its decomposition into prime powers as
$$n=\prod_{i =1}^m q_i^{a_i},$$ where the $q_i$'s are prime numbers, indexed by increasing order, and $a_i$'s are positive integers.
Denote by $p_i$ the $ i^{th }$ prime number, and for any integer $n,$ let
$$\bar{n}=\prod_{i =1}^m p_i^{a_i}.$$
Note that, by definition, for each $i=1,2,\cdots,m$ we have $q_i\ge p_i,$ and that, therefore, $n\ge \bar{n}.$
With this notation observe that 
$$\sigma(\bar{n})=\prod_{i= 1}^m\frac{p_i^{a_i+1}-1}{p_i-1}=\bar{n}\prod_{i= 1}^m\frac{p_i-p_i^{-a_i}}{p_i-1}.$$
In particular $$d(\bar{n})\le \prod_{i= 1}^m\frac{p_i}{p_i-1}\le 2^m,$$
and, likewise, $\frac{\sigma(n)}{n}\le 2.$ Thus, if $m$ is bounded
and $n \rightarrow \infty,$ we see that $d(n) \rightarrow \infty.$ We can thus assume when considering
$\liminf_{n \rightarrow \infty} d(n)$
that $m \rightarrow \infty.$
We prepare for the proof by a series of Lemmas.

{\lem For any integer $n\ge 1,$ we have $d(n)\ge d(\bar{n}).$}

\begin{proof}
Let $d(n)=f_1(n)-f_2(n),$ with $f_1(n)=e^\gamma \log \log n,$ and $f_2(n)=\frac{\sigma(n)}{n}.$
The monotonicity of the $\log$ and $n\ge \bar{n}$ yields $f_1(n)\ge f_1(\bar{n}).$
Write $f_2(n)=\prod_{i =1}^mg(a_i,q_i),$ where $g(a,x)=\frac{x-x^{-a}}{x-1}.$ Writing
$$g(a,x)=\frac{1+x+\cdots+x^a}{x^a}=\sum_{i=0}^a \frac{1}{x^i} ,$$
we see that, for fixed $a,$ the function $ x \mapsto g(a,x)$ is
nonincreasing in $x.$ This implies that $g(a_i,q_i)\le g(a_i,p_i)$ for each $i=1,2,\cdots,m$ and, therefore,
multiplying $m$ inequalities between nonegative numbers, that $f_2(n)\le f_2(\bar{n}).$ The result follows then by $d(n)=f_1(n)-f_2(n).$

\end{proof}

{\lem For any $n$ large enough we have $\frac{\sigma(\bar{n})}{\bar{n}}<e^\gamma \log p_m (1+\frac{1}{\log ^2 p_m}).$}

\begin{proof} Note that, with the notation of the proof of Lemma 1, we have $g(a,x) \le \frac{x}{x-1},$ for $x\ge 2$ and $a\ge 1,$ and, therefore
$$f_2(n)=\prod_{i =1}^mg(a_i,q_i)\le \prod_{i =1}^m\frac{p_i}{p_i-1}.$$

 The result follows then by \cite[Th. 8, (39)]{RS}.

\end{proof}

Recall Chebyshev summatory function $\vartheta(x)=\sum_{p\le x}\log(p).$ 

{\lem For all $n\ge 1,$ we have $\log \bar{n} \ge \vartheta(p_m).$}

\begin{proof}
 By definition $$ \log \bar{n} =\sum_{i= 1}^ma_i\log p_i\ge \sum_{i= 1}^m\log p_i=\vartheta(p_m).$$
\end{proof}

A classical result, related to the Prime Number Theorem, is

{\lem For large $x$ we have $\vartheta(x)=x+O(\frac{x}{\log x}).$}

\begin{proof}
An effective version is in \cite[Th. 4]{RS}. See for instance \cite[Th 4.7]{EM} for a sharper error term in $O(x\exp(-\frac{\sqrt{\log x}}{15})).$
\end{proof}

We are now ready for the proof  of Theorem 1. 
\begin{proof}

By Lemma 1 $d(n)\ge d(\bar{n}).$ %We only need to show that $\lim \frac{D(\bar{n})}{\bar{n}} = 0.$
By Lemma 2 we have 
\begin{equation}\label{lem2}
 -\frac{\sigma(\bar{n})}{\bar{n}}>-e^\gamma \log p_m (1+\frac{1}{\log ^2 p_m}).
\end{equation}

By Lemma 3 and 4 we have
\begin{equation}
 e^\gamma \log \log \bar{n}\ge e^\gamma \log \vartheta (p_m)=e^\gamma \log\left(p_m+O(\frac{p_m}{\log p_m})\right)=
 \end{equation}
 
\begin{equation}\label{lem34}
e^\gamma \left( \log p_m+\log(1+O(\frac{1}{\log p_m}))\right)=e^\gamma\left(\log(p_m)+O(\frac{1}{\log p_m})\right),
\end{equation}
where the last equality results from $\log (1+u) \sim u$ for $u \rightarrow 0.$
Adding up inequations \ref{lem2} and \ref{lem34}, after cancellation of the terms in $\log p_m,$ we obtain the inequality
$$d(\bar{n})=e^\gamma \log \log \bar{n}-\frac{\sigma(\bar{n})}{\bar{n}}\ge  O(\frac{1}{\log p_m})-\frac{e^\gamma}{\log p_m},$$
the right hand side of which goes to zero for large $n.$
\end{proof}
%%%%%%%%%%%%%%%%%%%%%%%%%%%%%%%%%%%%
\section{Proof of Theorem 3}
Recall the standard notation for oscillation theorems \cite[p. 194]{EM}.
If $f,g$ are two real valued functions of a real variable $x,$ with $g>0,$ then we write
\begin{itemize}
 \item $f(x)=\Omega_+(g(x)),$ if $\limsup_{x \rightarrow \infty}f(x)/g(x)>0$
 \item $f(x)=\Omega_-(g(x)),$ if $\liminf_{x \rightarrow \infty} f(x)/g(x)<0$
 \item $f(x)=\Omega_{\pm}(g(x)),$ if both $f(x)=\Omega_+(g(x)),$ and $f(x)=\Omega_-(g(x))$ hold
\end{itemize}

We refer the reader to \cite{R} for the definition of Colossally Abundant (CA) numbers.
By \cite[Proposition,\S 4]{R} if RH is false then, for CA numbers we have
$$ D(n)=\Omega_{\pm}( \frac{n \log \log n}{ (\log n)^b} ),$$
for some $b \in (0,1).$ This would imply, using the infinitude of CA numbers \cite{R}, that
$\liminf_{n \rightarrow \infty} D(n)= -\infty.$ 

If RH holds then by \cite[p.25]{Ra} the sequence $\frac{D(n)\sqrt{\log n}}{{n}}$ admits upper and lower limits for $n$ CA that are finite and $>0.$ Thus there are reals $>0$ say $A,\,B$ such that
$$A \frac{n}{\log n} \le D(n)\le B \frac{n}{\log n},$$ when $n$ is CA.
Therefore $\lim_{n \rightarrow \infty}D(n)=\infty.$ 
%%%%%%%%%%%%%%%%%%%%%%%%%%%%%%%%%%%%%%%%%%%%%%%%%%%%%%%%%%%%
\section{Conclusion and open problems}
In this note we have studied the quantity $D(n)$ which is the difference between the two handsides of Robin inequality, and its normalization $d(n)=\frac{D(n)}{n}.$
While the asymptotic behavior of $d(n)$ can be determined unconditionally (Theorem 1), that of $D(n)$ depends crucially on the truth of RH (Theorem 3).
It would be desirable to extend Theorem 3 to integers that are not CA.
It seems impossible to use Theorem 1 and Theorem 3 together to prove that RH holds. For instance, one cannot rule out the case that $D(n)$ behaves like
$-\sqrt{n}$ when $n\rightarrow \infty,$ which would not contradict the fact that $\liminf_{n \rightarrow \infty} d(n)=0.$

{\bf Acknowledgement:} Both authors are grateful to Minjia Shi for putting them in touch. They thank Gourab Batthacharya, Florian Luca, 
Jean-Louis Nicolas, Pieter Moree for helpful discussions.

\end{document}